\input xy
\xyoption{all}

\documentclass{jhrs}

\usepackage{amsmath}
\usepackage{xy}
\newtheorem {definition}[subsection]{Definition}
\newtheorem {lemma}[subsection]{Lemma}
\newtheorem {theorem}[subsection]{Theorem}

\begin{document}
\title{Separable morphisms of simplicial sets}          % Title of Document

\author{Dimitri Chikhladze}             % First Author
\email{d.chikhladze@gmail.com}       % First Author's email
\address{Department of Mathematics\\ % First Author's postal address
         Tbilisi State University\\
         2 University St., 0143 Tbilisi, Georgia}

% AMS 2000 Mathematics Subject Classification
\classification{55U10, 18A20, 18A40, 18G55, 57M10.}

% Keywords of the article
\keywords{simplicial set, separable morphism, covering map.}

% Abstract comes before maketitle
\begin{abstract}
We show that the class of separable morphisms in the sense of G.
Janelidze and W. Tholen in the case of Galois structure of second
order coverings of simplicial sets due to R. Brown and G. Janelidze
coincides with the class of covering maps of simplicial sets.

\end{abstract}

\received{}   % receive date (for example: 11 October 1999)
\revised{}    % receive date
\published{}  % publish date
\submitted{George Janelidze}  % Name of Journal's Editor, who submitted Article

\volumeyear{2006} % Volume Year
\volumenumber{1}  % Volume Number
\issuenumber{1}   % Issue Number

\startpage{1}     % PageNumber of first page

\maketitle

\section*{Introduction}
Separable morphisms were introduced in \cite{8} by A. Carboni and G.
Janelidze for lextensive categories. In the way of \cite{8} one can
consider separable morphism in a lextensive category $Fam(A)$, the
category of families of objects in a category $A$. What we call
$\Gamma_1$ below can be seen as a special case of this situation, a
characterization of separable morphisms for which is given by
Theorem \ref{teorema}.

G. Janelidze and W. Tholen \cite{4} defined separable morphisms in a
category $C$ for a given pointed endofunctor of $C$. Given an
adjunction $I, H : X \leftrightarrows C$ one can consider separable
morphisms with respect to the induced monad. Then, in a special
case, an appropriate adjunction $I, H : Sets \leftrightarrows
Fam(A)$ between $Fam(A)$ and the category of sets gives the same
notion of separable morphism in a lextensive category $Fam(A)$ as
\cite{8}.

Definition \ref{ganmarteba} in this paper is essentially of
\cite{4}, the difference being that in place of an adjunction we
consider a Galois structure, which together with a pair of adjoint
functors $I, H : X \leftrightarrows C$ consists of specified classes
of morphisms $F$ and $F^\prime$, called fibrations, in $C$ and $X$
respectively (see G. Janelidze \cite{6}, the earlier reference is G.
Janelidze \cite{5}), and we require separable morphisms to be
fibrations.

Our purpose is to describe the class of separable morphisms for the
Galois structure introduced by R. Brown and G. Janelidze in \cite{1}
($\Gamma_2$ below). Theorem \ref{teorema2} states that for this
Galois structure separable morphisms are exactly the Kan fibrations
which are covering maps of simplicity sets.

\section{Separable morphisms}\label{section}
In this section $C$ is a finitely complete category. Let $X$ be a
full reflective subcategory of $C$ with the inclusion $H: X
\rightarrow  C$. Suppose the reflection $I: C  \rightarrow X$ and
its unit  $\eta : 1  \rightarrow HI$ are chosen in a such way that
the counit is an identity $IH = 1$. Let $F$ and $F^\prime$ be
pullback stable classes of morphisms in $C$ and $X$ respectively,
such that $I(F) \subset F^\prime$  and $H(F^\prime) \subset F$.

Throughout this paper by Galois structure we mean the data $\Gamma =
(C, X, H, I,\eta,\\  F, F^\prime)$, with $C, X, H, I, \eta, F$ and
$F^\prime$  as above. Morphisms from $F$ and $F^\prime$ will be
called fibrations.

We will say that a morphism $f : A \rightarrow  B$ of the category
$C$ is a trivial covering or a cartesian morphism with respect to
the Galois structure  $\Gamma$  if it is a fibration, and the square

\medskip
$\xymatrix{A \ar[d]_f
 \ar[r]^{\eta_A} & HIA \ar[d]^{HIf} \\
B \ar[r]^{\eta_B} & HIB}$
\,\,\,\,\,\,\,\,\,\,\,\,\,\,\,\,\,\,\,\,\,\,\,\,\,\,\,\,\,\,\,\,\,\,\,\,\,\,\,\,\,\,\,\,\,\,\,\,\,\,\,\,\,\,\,\,\,\,\,\,\,\,\,\,\,\,\,\,\,\,\,\,\,\,\,\,\,\,\,\,\,\,\,\,\,\,\,\,\,\,\,\,\,\,\,\,\,\,\,\,\,\,\,\,\,\,\,\,\,\,\,\,\,\,\,\,\,\,\,\,\,\,\,\,\,\,\,\,\,\,\,\,\,\,\,\,\,\,\,\,\,\,\,\,\,\,\,(*)
\medskip

\noindent is a pullback (see \cite{3}, 3.1).

\begin{definition}\label{ganmarteba} A fibration $h : A  \rightarrow B$ from $F$ is called
a separable morphism with respect to the Galois structure $\Gamma =
(C, X, H, I, \eta , F, F^\prime )$  if the diagonal   $ \Delta =
\langle1_A, 1_A \rangle : A  \rightarrow A \times_B A$ is a
cartesian morphism with respect to $\Gamma$. \end{definition}

Consider an example. Suppose $A$ is a category with a terminal
object $t$. Let $Fam(A)$ be the category of families of objects in
$A$. A morphism $(m, \lambda) : (A_\lambda)_{\lambda \in \Lambda}
\rightarrow ({A^\prime}_{\lambda^\prime})_{\lambda^\prime \in
\Lambda^\prime}$ consists of a map $m : \Lambda \rightarrow
\Lambda^\prime$ and morphisms $\alpha_\lambda : A_\lambda
\rightarrow {A^\prime}_{f(\lambda)}$ for all $\lambda \in \Lambda$.
Let:

$C = Fam(A)$;

$X = Sets$ is the category of sets;

$H : X \rightarrow (A_x)_{x \in X}$ where $A_x = t$ for all $x$.

$I : (A_\lambda)_{\lambda \in \Lambda} \rightarrow \Lambda$; $I$ is
a left adjoint to $H$ with the obvious unit  $\eta$;

$F$ and $F^\prime$ are the classes of all morphisms in $C$ and $X$
respectively.

It is straightforward that a morphism $(m, \lambda) :
(A_\lambda)_{\lambda \in \Lambda} \rightarrow
({A^\prime}_{\lambda^\prime})_{\lambda^\prime \in \Lambda^\prime}$
is cartesian with respect to this Galois structure if and only if
all $\alpha_\lambda : A_\lambda \rightarrow {A^\prime}_{f(\lambda)}$
are isomorphisms.

\section{Separable morphisms for simplicial sets}
We will consider two Galois structures; first  $\Gamma_1$, in which:

$C = Sets^{\Delta^{op}}$ is the category of simplicial sets;

$X = Sets$ is the category of sets;

$H : X \rightarrow  C$ is the canonical inclusion;

$I : C  \rightarrow X$ is the functor sending a simplicial set to
the set of its connected components; it is a left adjoint to the
inclusion $H$ with the obvious unit  $\eta$;

$F$ and $F^\prime$ are the classes of all morphisms in $C$ and $X$
respectively.

Note that this structure is a special case of the example from the
section \ref{section}, when $A$ is the category of connected
simplicial sets. Then, it is easy to see that cartesian morphisms
with respect to $\Gamma_1$ are exactly the trivial coverings of
simplicial sets in the usual sense.

\begin{theorem}\label{teorema} A morphism of simplicial sets is a separable
morphism with respect to the Galois structure $\Gamma_1$ if and only
if for each commutative diagram

\medskip
$\xymatrix{\Delta[0] \ar[d]
 \ar[r] & A \ar[d] \\
\Delta[n] \ar@<0.7ex>[ur] \ar[ur] \ar[r] & B}$
\,\,\,\,\,\,\,\,\,\,\,\,\,\,\,\,\,\,\,\,\,\,\,\,\,\,\,\,\,\,\,\,\,\,\,\,\,\,\,\,\,\,\,\,\,\,\,\,\,\,\,\,\,\,\,\,\,\,\,\,\,\,\,\,\,\,\,\,\,\,\,\,\,\,\,\,\,\,\,\,\,\,\,\,\,\,\,\,\,\,\,\,\,\,\,\,\,\,\,\,\,\,\,\,\,\,\,\,\,\,\,\,\,\,\,\,\,\,\,\,\,\,\,\,\,\,\,\,\,\,\,\,\,\,\,\,\,\,\,\,\,\,\,\,\,\,\,(**)

\medskip

\noindent the diagonal morphisms are equal.
\end{theorem}

\begin{proof}
Observe that an injection $m : A \rightarrow B$ is cartesian with
respect to $\Gamma_1$ (i.e. trivial covering of simplicial sets) if
and only if each connected component of $B$ either is contained in
the image of $A$ under $m$, or has no intersection with it. We will
use this observation for $ \Delta : A \rightarrow A \times_B A$
(here and further $A \times_B A$ denotes the pullback of a typical
$h : A  \rightarrow B$ by itself; $\Delta$ is the diagonal $ \Delta
= \langle1_A, 1_A \rangle : A  \rightarrow A \times_B A$). In this
case, let $D$ denote the image of $A$ under $\Delta$; it consists of
all simplices of the form $(x, x)$ in $A \times_B A$.

Let in a diagram (**) $h$ be a separable morphism with respect to
$\Gamma_1$, what is to say that $\Delta$ is cartesian with respect
to $\Gamma_1$. If $x_1$ and $x_2$ are the n-simplices of $A$
corresponding to the diagonal morphisms of (**), then $h(x_1) =
h(x_2)$, and the pair $(x_1 ,x_2)$ is an n-simplex of $A \times_B
A$. The upper horizontal morphism in the diagram (**) gives a vertex
$a$ of $A$, and by commutativity of (**) $(a, a)$ is a vertex of
$(x_1, x_2)$. So, the connected component containing the simplex
$(x_1, x_2)$ has an intersection with $D$, but then by separability
of $\Delta$ it is completely contained in $D$. It follows that $x_1
= x_2$.

Now suppose that for the fixed $h : A \rightarrow B$ in each
commutative diagram (**) the diagonals are equal. We will prove that
$ \Delta : A \rightarrow A \times_B A$ is cartesian with respect to
$\Gamma_1$.

Let $K$ be a connected component of $A \times_B A$ the intersection
of which with $D$ is not empty. $K$ contains a vertex of the form
$(a, a)$, $a \in A[0]$. If $(d_1, d_2)$ ($d_1, d_2 \in A[0]$) is
another vertex of $K$, then there is a path from $(a, a)$ to $(d_1,
d_2)$ in $A \times_B A$.

Here we observe, if $(e_1, e_2$) is a 1-simplex in $A \times_B A$
``connecting" $(a^\prime, a^\prime)$ with $({d_1}^\prime,
{d_2}^\prime)$, then there is a commutative diagram

\medskip
$\xymatrix{\Delta[0] \ar[d]
 \ar[r] & A \ar[d] \\
\Delta[1] \ar@<0.7ex>[ur] \ar[ur] \ar[r] & B}$

\medskip

\noindent wherein the upper horizontal morphism is determined by
$a^\prime$ and the diagonals by $e_1$ and $e_2$. Since these
diagonals are equal $e_1 = e_2$, yielding also ${d_1}^\prime =
{d_2}^\prime$. Applying this argument consecutively to 1-simplices
in the chain connecting $(a, a)$ with $(d_1, d_2)$ we get $d_1 =
d_2$. Thus, each vertex of $K$ is in $D$.

Suppose now a pair $(x_1, x_2)$, with $x_1, x_2 \in A[n], n \geq 1$,
is an n-simplex of $K$. Since the vertices of $K$ are in $D$ we can
construct a commutative diagram of the form (**), wherein the
diagonals correspond to $x_1$ and $x_2$. These diagonals are equal,
so $x_1 = x_2$. We have proved that $K$ is completely contained in
$D$ showing that $\Delta$ is cartesian with respect to $\Gamma_1$.
\end{proof}

Consider now a Galois structure  $\Gamma_2$, introduced in \cite{1}:

$C = Sets^{\Delta^{op}}$ is the category of simplicial sets;

$X$ is the category of groupoids;

$H :  X \rightarrow C$ is the canonical inclusion called the nerve
functor;

$I = \pi_1 : C \rightarrow  X$ is the fundamental groupoid functor,
which is a left adjoint to $H$ with the obvious unit $\eta$;

$F$ and $F^\prime$  are the classes of Kan fibrations (\cite{3}, p.
65) in $C$ and $X$ respectively.

\begin{definition}\label{damf} A morphism of simplicial sets $h : A \rightarrow  B$ is called a covering
map if there exists a surjection $p : E  \rightarrow B$ and the
pullback of $h$ along $p$ is a trivial covering of simplicial sets.
\end{definition}

This definition is a special case of Definition 4.1 in \cite{6} for
the Galois structure $\Gamma_1$. However, as mentioned in
A.3.9(iii), \cite{7} by F. Borceux and G. Janelidze, it is
equivalent to the definition of covering map of simplicial sets
given in  \cite{3} by P. Gabriel and M. Zisman, where a morphism of
simplicial sets $h : A \rightarrow B$ is said to be a covering map
if and only if for each commutative square

\medskip

$\vcenter{\xymatrix{\Delta[0] \ar[d]_i
 \ar[r]^v & A \ar[d]^h \\
\Delta[n] \ar[r]^u & B}}$
\,\,\,\,\,\,\,\,\,\,\,\,\,\,\,\,\,\,\,\,\,\,\,\,\,\,\,\,\,\,\,\,\,\,\,\,\,\,\,\,\,\,\,\,\,\,\,\,\,\,\,\,\,\,\,\,\,\,\,\,\,\,\,\,\,\,\,\,\,\,\,\,\,\,\,\,\,\,\,\,\,\,\,\,\,\,\,\,\,\,\,\,\,\,\,\,\,\,\,\,\,\,\,\,\,\,\,\,\,\,\,\,\,\,\,\,\,\,\,\,\,\,\,\,\,\,\,\,\,\,\,\,\,\,\,\,\,\,\,\,\,\,\,\,\,\,\,(***)

\medskip

\noindent there exists a unique morphism $x :  \Delta[n] \rightarrow
A$ with $hx = u$ and $xi = v$; $x$ satisfying these equations is
called a diagonal fill-in of the square (***).

\begin{lemma}\label{lemma}
An injection of simplicial sets is cartesian with respect to
$\Gamma_2$ if and only if it is cartesian with respect to
$\Gamma_1$.
\end{lemma}

\begin{proof}
A trivial covering of simplicial sets is a Kan fibration. Also, it
is not difficult to see that for the Galois structure $\Gamma_2$, a
square (*) is a pullback if $h$ is a trivial covering of simplicial
sets. The ``if" part of the Lemma is proved.

Conversely, if $m : A \rightarrow B$ is cartesian with respect to
$\Gamma_2$, then, by definition, it is a Kan fibration. In
particular, the inverse image of each connected component $K$ of $B$
is either empty or surjectively mapped on $K$ by $m$. If in addition
to this $m$ is injective, then clearly it is a trivial covering of
simplicial sets.
\end{proof}

\begin{theorem}\label{teorema2}
A Kan fibration of simplicial sets is separable with respect to
$\Gamma_2$ if and only if it is a covering map.

\end{theorem}

\begin{proof}
By Lemma \ref{lemma} the diagonal $ \Delta : A \rightarrow A
\times_B A$ is cartesian with respect to $\Gamma_2$ if and only if
it is cartesian with respect to $\Gamma_1$. Then, Theorem
\ref{teorema2} will follow from Theorem \ref{teorema} if we note
that in a commutative diagram (***) where $h$ is a Kan fibration a
diagonal fill-in exists.
\end{proof}

\end{document}